\documentclass[square, comma]{article}
\usepackage[utf8]{inputenc}
\usepackage{amsmath}
\usepackage{amsthm}
\usepackage{cite}
\usepackage{hyperref}

\title{An Elementary Proof of the Prime Number Theorem based on Möbius Function }
\author{Junda Pan }
\date{\emph{email: jundapan@foxmail.com}}
\begin{document}

\maketitle

\begin{abstract}
    Let $\mu(n)$ denote the Möbius function, define $M(x)= \sum_{n\leq x}^{}\mu (n)$. The main result of this paper is to prove that 
    \begin{equation*}
        \displaystyle\lim_{x \to +\infty}\frac{M(x)}{x}=0
    \end{equation*}
    which is equivalent to the prime number theorem. We also use Selberg's asymptotic formula, but the treatments of key parts are different from several classical proofs.
\end{abstract}

{\centering\section{Introduction}}
The prime number theorem, one of the most well-known and significant theorems in number theory, was first proved by Hadamard and de la Vallée Poussin in 1896 using analytic methods. The elementary proof was discovered in 1949 by A. Selberg \cite{ref1} and P. Erdős \cite{ref2}. Their proof makes no use of complex analysis nor of the Riemann zeta function but is intricate. In addition, there are more classic elementary proofs such as \cite{ref3, ref4, ref5, ref6}.\par
Let 
\begin{equation*}
\mathcal{F}(x):=\sum_{n\leq x}^{}\mu(n)\mathrm{log}\ \frac{x}{n}=\int_{1}^{x}\frac{M(y)}{y}\mathrm{d}y\ ,\ 
\mathcal{H}(x):=\frac{\mathcal{F}\left ( e^{\sqrt{x}} \right )}{e^{\sqrt{x}}}.
\end{equation*}
Then we can assume that exist positive constant $\ell,L,\alpha $ such that
\\
\begin{equation*}
\varlimsup\limits_{}|\mathcal{H}(x)|=\ell\ ,\ 
\varlimsup\limits_{}\frac{1}{x}\int_{0}^{x}|\mathcal{H}(y)|\mathrm{d}y=L,
\end{equation*}
\begin{equation}
|\mathcal{H}(x)|\leq\alpha + o(1),
\end{equation}
by noticing that $|\mathcal{H}(x)|\leq1$.
To derive
\begin{equation}
\displaystyle\lim_{x \to +\infty}\frac{M(x)}{x}=0,
\end{equation}
we shall prove $\ell=0$ in two cases, which can be proved to be equivalent to (2). In one case, we verified that $\ell=0$ directly by some simple skills in calculus; in another case, we got an inequality
\begin{equation*}
\ell\leq\frac{1}{1+\lambda }\alpha ,
\end{equation*}
where $\lambda\in \left ( 0,1 \right ) $ is a constant.
Then we can easily obtain $\ell=0$ through iteration.\par
In this paper, we omit many computational details because we believe that such processing can enable readers to grasp the framework and ideas of the proof faster. In our proof, Lemma 2 is based on N. Levinson \cite[part 3]{ref4}. The whole method, of course, still belongs to Selberg in essence.  
\\

{\centering\section{Preliminaries}}
We list some well-known or trivial results without proof.\\
\\
\textbf{Proposition 1} (Tatuzawa, Iseki). \emph{Let F be a real- or complex-valued function defined on} $(1,\infty)$, $\Lambda$ \emph{be Mangoldt function,} \emph{and let }\\
\begin{equation*}
G(x)=\mathrm{log}\ x\sum_{n\leq x}^{}F\left ( \frac{x}{n} \right ).
\end{equation*}
\emph{Then}
\begin{equation*}
F(x)\mathrm{log}\ x+\sum_{n\leq x}^{}F\left ( \frac{x}{n} \right )\Lambda (n)=\sum_{d\leq x}^{}\mu (d)G\left ( \frac{x}{d} \right ).
\end{equation*}
\\
\textbf{Proposition 2 }(Selberg). \emph{Let $\Lambda_2:=\mu \ast (\mathrm{log})^2=\Lambda \ast \Lambda  +\Lambda \mathrm{log}$, thus}

\begin{equation}
\sum_{n\leq x}^{}\Lambda _2(n)=2x\ \mathrm{log}\ x+\mathcal{O}(x).
\end{equation}
Defining that $\Theta=(\Lambda \ast \Lambda ) /\mathrm{log}$, then (3) can be written as

\begin{equation}
\sum_{n\leq x}^{}(\Lambda (n)+\Theta (n))=2x+\mathcal{O}(\frac{x}{\mathrm{log}\ x})
\end{equation}

\noindent by partial summation.\\


\noindent \textbf{Proposition 3.} \emph{We have}

\begin{equation}
\sum_{n\leq x}^{}\mathcal{F}\left ( \frac{x}{n} \right )=\sum_{n\leq x}^{}\mu (n)\biggl\lfloor\frac{x}{n}\biggr\rfloor\mathrm{log}\ \frac{x}{n}=\mathcal{O}(\mathrm{log}\ x).
\end{equation}

Then
\begin{equation*}
\int_{1}^{x}\mathcal{F}\left ( \frac{x}{t} \right )\mathrm{d}t=\mathcal{O}(x),
\end{equation*}
since
\begin{equation*}
\mathcal{F}\left ( \frac{x}{n} \right )=\int_{n}^{n+1}\mathcal{F}\left ( \frac{x}{t} \right )\mathrm{d}t+\mathcal{O}\left ( \sum_{x/(n+1)<k<x/n}^{}1 \right ).
\end{equation*}
It gives
\begin{equation}
\left|\int_{x_1}^{x_2}{}\mathcal{H}(y)\mathrm{d}y \right|\leq M,
\end{equation}
where $M$ is a constant, $x_1$, $x_2$ are arbitrary non-negative numbers.\\

\noindent \textbf{Proposition 4.} \emph{Assume $f(x)\geq0,\ f(a)=0,\ |f'(x)|\leq m$. Then}

\begin{equation*}
\int_{a}^{b}f(x)\mathrm{d}x\leq \frac{1}{2}m(b-a)^2.
\end{equation*}
\\
{\centering\section{Lemmas}}
\noindent \textbf{Lemma 1.} \emph{For $x>0$ we have}
\begin{equation}
|\mathcal{H}(x)|\leq\frac{1}{x}\int_{0}^{x}|\mathcal{H}(y)|\mathrm{d}y+\mathcal{O}\left ( \frac{1}{\sqrt{x}} \right ).
\end{equation}
\begin{proof}
We start with Proposition 1. Let $F=\mathcal{F}$, (5) gives
\begin{equation}
\mathcal{F}(x)\mathrm{log}\ x+\sum_{n\leq x}^{}\mathcal{F}\left ( \frac{x}{n} \right )\Lambda (n)=\mathcal{O}(x)
\end{equation}

\noindent by noticing that
\begin{equation}
\sum_{n\leq x}^{}\left ( \mathrm{log}\ \frac{x}{n} \right )^2=2x+\mathcal{O}((\mathrm{log}\ x)^2).
\end{equation}

\noindent Replacing $x$ with $x/n$, then (8) becomes
\begin{equation}
\mathcal{F}\left ( \frac{x}{n} \right )\mathrm{log}\ \frac{x}{n}+\sum_{m\leq x/n}^{}\mathcal{F}\left ( \frac{x}{mn} \right )\Lambda (m)=\mathcal{O}\left ( \frac{x}{n} \right ).
\end{equation}
Noticing that
\begin{equation*}
\mathcal{O}\left ( x \right )\mathrm{log}\ x-\sum_{n\leq x}^{}\Lambda (n)\mathcal{O}\left ( \frac{x}{n} \right )=\mathcal{O}(x\ \mathrm{log}\ x),
\end{equation*}
since
\begin{equation*}
\sum_{n\leq x}^{}\frac{\Lambda (n)}{n}=\mathrm{log}\ x+\mathcal{O}(1).
\end{equation*}
So we combine (8), (10)
\begin{equation*}
\left\{\mathcal{F}(x)\mathrm{log}\ x+\sum_{n\leq x}^{}\mathcal{F}\left ( \frac{x}{n} \right )\Lambda (n) \right\}\mathrm{log}\ x-
\end{equation*}
\begin{equation*}
\sum_{n\leq x}^{}\Lambda (n)\left\{\mathcal{F}\left ( \frac{x}{n} \right )\mathrm{log}\ \frac{x}{n}+\sum_{m\leq x/n}^{}\mathcal{F}\left ( \frac{x}{mn} \right )\Lambda (m) \right\}=\mathcal{O}(x\ \mathrm{log}\ x).
\end{equation*}
It follows that
\begin{equation}
\mathcal{F}(x)\left ( \mathrm{log}\ x \right )^2=\sum_{n\leq x}^{}\Lambda ^-_2\mathcal{F}\left ( \frac{x}{n} \right )+\mathcal{O}(x\ \mathrm{log}\ x),
\end{equation}
or
\begin{equation}
|\mathcal{F}(x)|\left ( \mathrm{log}\ x \right )^2\leq \sum_{n\leq x}^{}\Lambda _2\left|\mathcal{F}\left ( \frac{x}{n} \right ) \right|+\mathcal{O}(x\ \mathrm{log}\ x),
\end{equation}
where
\begin{equation*}
\Lambda^-_2:=\Lambda\ast \Lambda -\Lambda \mathrm{log}.
\end{equation*}
From Proposition 2, we consider
\begin{equation*}
\sum_{n\leq x}^{}\left\{ \Lambda ^-_2(n)-2\Lambda (n)\mathrm{log}\ \frac{x}{n}\right\}=\mathcal{O}(x\ \mathrm{log}\ x).
\end{equation*}
It follows that
\begin{equation*}
\Lambda ^-_2(n)-2\Lambda (n)\mathrm{log}\ \frac{x}{n}=\mathcal{O}(\mathrm{log}\ n),
\end{equation*}
or
\begin{equation}
2\Lambda (n)\mathrm{log}\ \frac{x}{n}=\left| \Lambda ^-_2(n)\right|+\mathcal{O}(\mathrm{log}\ n).
\end{equation}
Similarly, we have
\begin{equation}
2\ \mathrm{log}\ n= \Lambda _2(n)+\mathcal{O}(1).
\end{equation}
Then (5), (12), (14) yield
\begin{equation}
|\mathcal{F}(x)|\left ( \mathrm{log}\ x \right )^2\leq 2\sum_{n\leq x}^{} \left|\mathcal{F}\left ( \frac{x}{n} \right ) \right|\mathrm{log}\ n+\mathcal{O}(x\ \mathrm{log}\ x).
\end{equation}
Using (11), (13), (15) we obtain
\begin{equation}
|\mathcal{F}(x)|\left ( \mathrm{log}\ x \right )^2\leq 2\sum_{n\leq x}^{}\Lambda (n) \left|\mathcal{F}\left ( \frac{x}{n} \right ) \right|\mathrm{log}\ \frac{x}{n}+\mathcal{O}(x\ \mathrm{log}\ x).
\end{equation}
Similarly, we have
\begin{equation}
|\mathcal{F}(x)|\left ( \mathrm{log}\ x \right )^2\leq 2\sum_{n\leq x}^{}\Theta  (n) \left|\mathcal{F}\left ( \frac{x}{n} \right ) \right|\mathrm{log}\ \frac{x}{n}+\mathcal{O}(x\ \mathrm{log}\ x),
\end{equation}
just observe
\begin{equation*}
\Theta(n)+\Lambda (n)=\mathcal{O}(1),
\end{equation*}
which is a corollary of (4).\\
\\
Thus, by combining (16), (17) and using (4) again, we derive
\begin{equation}
|\mathcal{F}(x)|\left ( \mathrm{log}\ x \right )^2\leq 2\sum_{n\leq x}^{} \left|\mathcal{F}\left ( \frac{x}{n} \right ) \right|\mathrm{log}\ \frac{x}{n}+\mathcal{O}(x\ \mathrm{log}\ x),
\end{equation}
which follows by partial summation.
Similar to the operation in Proposition 3, (18) can be written as follows
\begin{equation*}
|\mathcal{F}(x)|\left ( \mathrm{log}\ x \right )^2\leq 2\int_{1}^{x}\left|\mathcal{F}\left ( \frac{x}{t} \right ) \right|\mathrm{log}\ \frac{x}{t}\mathrm{d}t+\mathcal{O}(x\ \mathrm{log}\ x).
\end{equation*}
Finally we have, on writing  $e^{\sqrt{x}}$ for $x$
\begin{equation*}
|\mathcal{H}(x)|\leq\frac{1}{x}\int_{0}^{x}|\mathcal{H}(y)|\mathrm{d}y+\mathcal{O}\left ( \frac{1}{\sqrt{x}} \right ).
\end{equation*}
This completes the proof.
\end{proof}

\noindent \textbf{Lemma 2.} \emph{Let $a,b$ be successive zeros of $\mathcal{H}(x)$, then there are positive constants $\kappa $, $\varepsilon$ such that}
\begin{equation}
\int_{a}^{b}|\mathcal{H}(x)|\mathrm{d}x\leq \alpha(b-a) (1-\kappa |\mathcal{H}(\xi)|  )+o(b-a-\varepsilon ),
\end{equation}
\emph{where $a< \xi < b$, $\alpha $ is defined in} Introduction.

\begin{proof} We take a point $a+\alpha /h$, where h is an arbitrary positive number such that
\begin{equation}
\int_{a}^{b}=\int_{a}^{a+\frac{\alpha }{h}}+\int_{a+\frac{\alpha }{h}}^{b}
\end{equation}
Noticing that $b-a\geq 1$, it can always have achieved because of the arbitrariness of $h$. May as well let $h=h_1$ satisfy (20).
Clearly, $\mathcal{H}'(x)$ is bounded, so we get from (1) and Proposition 4
\begin{equation*}
\int_{a}^{b}|\mathcal{H}(x)|\mathrm{d}x\leq \frac{m}{2}\left ( \frac{\alpha }{h} \right )^2+(b-a-\frac{\alpha }{h})\ \alpha +o(b-a-\frac{\alpha }{h})
\end{equation*}
\begin{equation}
=\alpha(b-a)\left ( 1+\frac{(m-2h)\alpha }{2h^2(b-a)} \right )+o(b-a-\frac{\alpha }{h}).
\end{equation}
\\
We choose $h=h_2$ such that $m-2h<0$.
By (6), existing $\xi $, $a<\xi<b$, such that
\begin{equation*}
\int_{a}^{b}|\mathcal{H}(x)|\mathrm{d}x=|\mathcal{H}(\xi )|\ (b-a)\leq M.
\end{equation*}
Hence
\begin{equation*}
\int_{a}^{b}|\mathcal{H}(x)|\mathrm{d}x\leq\alpha(b-a)\left ( 1-\frac{(2h-m)\alpha }{2Mh^2}|\mathcal{H}(\xi )| \right )+o(b-a-\frac{\alpha }{h}),
\end{equation*}
just assume $h=max\left\{h_1,h_2 \right\}$.\\
\\
Marking $\kappa=(2h-m)\alpha/\left (2Mh^2  \right ),\varepsilon =\alpha /h $ then we complete the proof.
\end{proof}

{\centering\section{Proof of (2)}}

\noindent \textbf{Theorem 1.} \emph{We have} $\ell=0$.
\\
\begin{proof} Now we discussing the zeros of $\mathcal{H}(x)$.
If $\mathcal{H}(x)$ has finite zeros, it implies $\ell=0$ directly since (6) and Lemma 1. So we assume that $\mathcal{H}(x)$ has infinite zeros. Let $\delta(x) $ denote the number of zeros not exceeding $x$, $x_i\ (i=1,2,\dots,\delta (x),\dots)$ denote the ith zero, then
\begin{equation*}
\frac{1}{x}\int_{0}^{x}|\mathcal{H}(y)|\mathrm{d}y=\frac{1}{x}\int_{0}^{x_1}+\int_{x_1}^{x_2}+\cdots+\int_{x_{\delta (x)-1}}^{x_{\delta (x)}}+\int_{x_{\delta (x)}}^{x}|\mathcal{H}(y)|\mathrm{d}y
\end{equation*}

\begin{equation}
\leq \frac{2M}{x}+\frac{1}{x}\int_{x_1}^{x_2}+\cdots+\int_{x_{\delta (x)-1}}^{x_{\delta (x)}}|\mathcal{H}(y)|\mathrm{d}y,
\end{equation}
and existing $\xi_i\  (i=1,2,\dots,\delta (x)-1)$ such that
\begin{equation*}
\int_{x_i}^{x_{i+1}}|\mathcal{H}(y)|\mathrm{d}y=\left ( x_{i+1}-x_i \right )|\mathcal{H}(\xi _i)|.
\end{equation*}
Thus, we divide it into two cases.\\
\par\emph{Case 1.} $\displaystyle \lim_{ n\to\infty }|\mathcal{H}(\xi _n)|=0$.\\
\\
Then when $x\to \infty$, we have
\begin{equation*}
\frac{1}{x_{\delta (x)}-x_{\delta (x)-1}}\int_{x_{\delta (x)-1}}^{x_{\delta (x)}}|\mathcal{H}(y)|\mathrm{d}y\to 0.
\end{equation*}
This gives
\begin{equation*}
\frac{\mathcal{G}\left (  x_{\delta (x)}\right )-\mathcal{G}\left ( x_{\delta (x)-1} \right )}{x_{\delta (x)}-x_{\delta (x)-1}}\to 0,
\end{equation*}
where 
\begin{equation*}
\mathcal{G}(x):=\int_{x_1}^{x}|\mathcal{H}(y)|\mathrm{d}y.
\end{equation*}
Hence
\begin{equation*}
\frac{\mathcal{G}(x_{\delta (x)})}{x_{\delta (x)}}\to 0.
\end{equation*}
It can be written as
\begin{equation*}
\frac{1}{x_{\delta (x)}}\int_{x_1}^{x_2}+\int_{x_2}^{x_3}+\cdots+\int_{x_{\delta (x)-1}}^{x_{\delta (x)}}|\mathcal{H}(y)|\mathrm{d}y\to 0.
\end{equation*}
By noticing that $x\geq x_{\delta (x)}$, we get from Lemma 1 and (22)
\begin{equation*}
0\leq \ell\leq L\leq 0
\end{equation*}
which implies $\ell=0$.\\
\par\emph{Case 2.} $\displaystyle \lim_{ n\to\infty }|\mathcal{H}(\xi _n)|>0$ or the limit does not exist.\\
\\
Then the minimum number exists in $\left\{\left|\mathcal{H}(\xi _i) \right|\mid i=1,2,\dots,x_{\delta (x)},\dots \right\}$. Let $\iota:=\underset{i}{min}\left|\mathcal{H}(\xi _i) \right| $ , we have
\begin{equation*}
\frac{1}{x}\int_{x_1}^{x_2}+\dots+\int_{x_{\delta (x)-1}}^{x_{\delta (x)}}|\mathcal{H}(y)|\mathrm{d}y\leq
\end{equation*}

\begin{equation*}
\frac{\alpha}{x}(x_{\delta (x)}-x_1) (1- \kappa\iota )+o\left\{x_{\delta (x)}-x_1-(\delta (x)-1)\varepsilon  \right\},
\end{equation*}
since Lemma 2.
It follows that
\begin{equation*}
\ell\leq L\leq \alpha -\lambda \alpha
\end{equation*}
where $\lambda=\kappa \iota  $, $0<\lambda<1$.
By $0\leq \ell\leq \lambda$, we get
\begin{equation}
0\leq \ell\leq \frac{1}{1+\lambda }\alpha .
\end{equation}
Repeating the above operation n times, then (23) becomes
\begin{equation*}
0\leq \ell\leq \frac{1}{\lambda_n }\alpha ,
\end{equation*}
where $\lambda_{n}=1+\lambda\lambda_{n-1}$, $\lambda _0=1$.
It implies $\ell=0$ when we let $n\to \infty$.\\
\par To sum up, we complete the proof.
\end{proof}

\noindent\textbf{Corollary 1.} \emph{Let $H(x):=e^{-\sqrt{x}}M\left ( e^{\sqrt{x}} \right )$, then}
\begin{equation*}
\displaystyle \lim_{x\to +\infty}H(x)=0.
\end{equation*}

\begin{proof} We have
\begin{equation}
|H(x)|\leq\frac{1}{x}\int_{0}^{x}|H(y)|\mathrm{d}y+o(1).
\end{equation}
The proof is completely contained in the proof of Lemma 1. Finally, we get from Theorem 1 and (24)
\begin{equation}
0\leq \varlimsup\limits_{}|H(x)|\leq \varlimsup\limits_{}|\mathcal{H}(x)|=0.
\end{equation}
It implies the corollary and (2) at once.
\end{proof}

\textsc{Final remark.} 
1. As one sees Lemma 1 is similar to \cite[(2.14)]{ref1}, but if we use the method in \cite{ref1}, we cannot get the required accuracy (but the method can be used for the proof of (24)). In reality, similar to the method in \cite{ref2}, if we only use (15), we can also derive (2) by dealing with a double integral.
2. (24) is necessary. We cannot get (2) from Theorem 1 directly since $M(x)$ is not monotonic.

\centering

\end{document}